\documentclass[12pt, a4paper]{article}
\usepackage{amssymb, amsmath, color}
\usepackage[latin1]{inputenc}
\usepackage{graphicx}
\usepackage{enumerate} 
\usepackage{color}

\newtheorem{theorem}{Theorem}
\newtheorem{proposition}[theorem]{Proposition}
\newtheorem{lemma}[theorem]{Lemma}
\newtheorem{remark}[theorem]{Remark}
\newtheorem{definition}[theorem]{Definition}
\newtheorem{corollary}[theorem]{Corollary}

\begin{document}

\begin{center}
{\Large \bf {Geometry of extended Bianchi-Cartan-Vranceanu spaces}}
\end{center}

\vspace{1cm}

\begin{center}
{\bf Angel Ferr\'andez$^{1*}$, Antonio M. Naveira$^2$ and Ana D. Tarrío$^3$}
\end{center}

\vspace{1cm}

\noindent {\small ${}^1$Departamento de Matem\'aticas, Universidad de Murcia,
Campus de Espinardo, 30100 Murcia, Spain.
E-mail address: {\ttfamily aferr@um.es}\\[1mm]
${}^2$Departamento de Matemáticas, Universidad de Valencia (Estudi General),
Campus de Burjassot, 46100 Burjassot, Spain.
E-mail address: {\ttfamily naveira@uv.es}\\[1mm]
${}^3$Departamento de Matemáticas, Universidade da Coruña, Campus A Zapateira,
15001 A Coruña, Spain.
E-mail address: {\ttfamily madorana@udc.es}}

\begin{abstract} The differential geometry of $3$-dimensional Bianchi, Cartan and Vranceanu ($BCV$) spaces is well known. We introduce the extended Bianchi, Cartan and Vranceanu ($EBCV$) spaces as a natural seven dimensional generalization of $BCV$ spaces and study some of their main geometric properties, such as the Levi-Civita connection, Ricci curvatures, Killing fields and geodesics.
\end{abstract}

\medskip



\def\thefootnote{}
\footnotetext{(*) Corresponding author A. Ferr\'andez.}


\section{Introduction}

Let us denote by $H^{2n+1}$ the $(2n + 1)$-dimensional complex Heisenberg group in
$\mathbb{R}^{2n+1}=\mathbb{C}^n \times \mathbb{R}$ with coordinates $(z,t)=(x, y, t) =
(x_1, y_1, \ldots , x_n, y_n, t)$ whose group law writes
$$(x,y,t)\cdot(x',y',t')=(z,t)\cdot(z',t')=(z+z',t+t'+\frac{1}{2}{\rm Im}\sum_{j=1}^n z_j{\bar z}'_j).$$
Set $g_0 = (z_0,t) \in H^{2n+1}$ and let $l_{g_0}g = g_0 g$ be the left translation by $g_0$. We can then easily see that the left invariant vector fields write down
$$X_{\alpha}=\frac{\partial}{\partial x_{\alpha}}+\frac{1}{2}y_{\alpha}\frac{\partial}{\partial t}, \quad
Y_{\alpha}=\frac{\partial}{\partial y_{\alpha}}-\frac{1}{2}x_{\alpha}\frac{\partial}{\partial t}, \quad T=\frac{\partial}{\partial t},$$ and $\{X_{\alpha}, Y_{\alpha}, {\alpha}=1, \ldots, n\}$ form an orthonormal basis of a distribution $D$ with respect to the sub-Riemannian metric ${\rm ds}^2 =\sum_{{\alpha}=1}^n (dx_{\alpha}^2 + dy_{\alpha}^2)$ and satisfy the following bracket relations:
$$[X_{\alpha},Y_{\alpha}] = T, \quad [X_{\alpha},T] = 0, \quad [Y_{\alpha},T] = 0, \quad {\alpha} = 1, \ldots, n.$$
Taking $n=1$ and $t \in \mathbb{S}^1$ we obtain the complex Heisenberg group, which is a manifold equipped with a contact structure.

The quaternionic Heisenberg group serves as a flat model of quaternionic contact manifolds. We can consider the following model $\mathbb{Q}^n \times {\rm Im}\mathbb{Q}$ with the group law
$$(q, p) = (q_1, p_1)\circ(q_2, p_2) = (q_1 + q_2, p_1 + p_2 + \frac{1}{2}{\rm Im}(p_1 \bar{p_2}),$$
where $q_1, q_2 \in \mathbb{Q}^n,\, p_1, p_2 \in {\rm Im}\mathbb{Q}$, and $\mathbb{Q}$ stands for the quaternionic field $\mathbb{Q} = \{ q = w + x I + y J + z K, (w, x, y, z) \in \mathbb{R}^4\}$ and ${\rm Im}\mathbb{Q} = \{ p = r i + s j + t k, (r, s, t) \in \mathbb{R}^3\}$ with the Pauli matrices

$$
I=\left(%
\begin{array}{cc}
i & 0 \\
0 & -i \\
\end{array}%
\right), \qquad
J=\left(%
\begin{array}{ll}
0 & 1 \\
-1 & 0 \\
\end{array}%
\right), \qquad
K=\left(%
\begin{array}{ll}
0 & i \\
i & 0 \\
\end{array}%
\right).
$$

Bearing in mind these elementary algebraic computations, it is easily understood the definitions of the Bianchi-Cartan-Vranceanu spaces ($BCV$ spaces for short) as well as their natural extensions, as we will do in next sections.

\section{The Bianchi-Cartan-Vranceanu (BCV) spaces (see \cite{Ca-Cha})}

It was Cartan (\cite{Ca}) who obtained the families of today known as $BCV$-spaces by classifying  three-dimensional Riemannian manifolds with four-dimensional isometry group. They also appeared in the work of L. Bianchi (\cite{Bi2, Bi3}), and G. Vranceanu (\cite{V}). These kind of spaces have been extensively studied and classified (see for instance \cite{Profir, VV}). In theoretical cosmology they are known as  Bianchi-Kantowski-Saks spaces, which are used to construct some homogeneous spacetimes (\cite{LL}).

For real numbers $m$ and $l$, consider the set
$$BCV(m,l)=\{(x, y, z) \in \mathbb{R}^3 \, : \, 1+m(x^2+y^2)>0\}$$ equipped with the metric
$${\rm ds}^2_{m,l} = \frac{dx^2+dy^2}{(1+m(x^2+y^2))^2} + \left(dr+ \frac{l}{2}\;\frac{x dy -ydx}{1+m(x^2+y^2} \right)^2.$$

Observe that this metric is obtained as a conformal deformation of the planar Euclidean metric by
adding the imaginary part of $z\, d \bar z$, for a complex number $z$.

The complete classification of $BCV$ spaces is as follows:

(i) If $m=l=0$, then $BCV(m,l)\cong \mathbb{R}^3$;

(ii) If $m=\frac{l}{4}$, then $BCV(m,l)\cong (\mathbb{S}^3(m)-\{\infty\})$;

(iii) If $m>0$ and $l=0$, then $BCV(m,l)\cong (\mathbb{S}^2(4m)-\{\infty\}) \times \mathbb{R}$;

(iv) If $m<0$ and $l=0$, then $BCV(m,l)\cong (\mathbb{H}^2(4m)-\{\infty\}) \times \mathbb{R}$;

(v) If $m>0$ and $l \neq 0$, then $BCV(m,l)\cong {\rm SU(2)}-\{\infty\}$;

(vi) If $m<0$ and $l \neq 0$, then $BCV(m,l)\cong \widetilde{{\rm SL}}(2,\mathbb{R})$;

(vii) If $m=0$ and $l \neq 0$, then $BCV(m,l)\cong {\rm Nil_3}$.

The following vector fields form an orthonormal frame of $BCV(m,l)$:

$$E_1=\left(1+m(x^2+y^2)\right)\partial_x-\frac{l}{2}y\partial_z,\qquad E_2=\left(1+m(x^2+y^2)\right)\partial_y+\frac{l}{2}x\partial_z,\qquad E_3=\partial_z.$$


Let $\mathcal{D}$ be the distribution generated by $\{E_1,E_2\}$, then the manifold $\left(BCV(m,l), \mathcal{D}, {\rm ds}^2_{m,l}\right)$ is an example of sub-riemannian geometry (see \cite{ Ca-Cha, Strichartz}) and the horizontal distribution is a 2-step breaking-generating distribution everywhere.


\section{Extended Bianchi-Cartan-Vranceanu spaces}\label{Sect2}

\subsection{Set up} Observe that letting $z=x+iy$, we see that ${\rm Im}(z\,d\bar z)=ydx-xdy$, which reminds us the map $\mathbb{C}\times \mathbb{C} \to \mathbb{R}\times \mathbb{C}$ given by $(z_1,z_2)\mapsto
(|z_1|^2-|z_2|^2, 2(z_1{\bar z}_2))$, that easily leads to the classical Hopf fibration $\mathbb{S}^1 \hookrightarrow \mathbb{S}^3 \to \mathbb{S}^2$, where coordinates in $\mathbb{S}^2$ are given by
$(|z_1|^2-|z_2|^2, 2{\rm Re}(z_1{\bar z}_2), 2{\rm Im}(z_1{\bar z}_2))$.

In the same line, using quaternions $\mathbb{H}$ instead of complex numbers, we get the  fibration $\mathbb{S}^3 \hookrightarrow \mathbb{S}^7 \to \mathbb{S}^4$. Quaternions are usually presented with the imaginary units $i, j, k$ in the form $q = x_0 + x_1i + x_2j + x_3k$, $x_0, x_1, x_2, x_3 \in \mathbb{R}$ with $i^2 = j^2 = k^2 = ijk =-1$. They can also be defined equivalently, using the complex numbers $c_1 = x_0 + x_1i$ and $c_2 = x_2 + x_3i$, in the form $q = c_1 + c_2 j$. Then for a point $(q_1 = \alpha+\beta j, q_2 =\gamma+\delta j) \in \mathbb{S}^7$, we get the following coordinate expressions $(|q_1|^2-|q_2|^2, 2{\rm Re}(\bar{\alpha}\gamma+\bar{\beta}\delta), 2{\rm Im}(\bar{\alpha}\gamma+\bar{\beta}\delta), 2{\rm Re}(\alpha \delta-\beta \gamma), 2{\rm Im}(\alpha \delta-\beta \gamma))$.

For any $q=w+xi+yj+zk \in \mathbb{H}$ we find that $qd\bar{q}=wdw+xdx+ydy+zdz+(xdw-wdx+zdy-ydz)i+(ydw-wdy+xdz-zdx)j+(zdw-wdz+ydx-xdy)k$. As the quaternionic contact group $\mathbb{H} \times {\rm Im}\mathbb{H}$, with coordinates $(w,x,y,z,r,s,t)$ can be equipped with the metric
\begin{equation*}
\begin{split}
{\rm ds}^2 &=
(dw^2+dx^2+dy^2+dz^2)
+\left(dr+\frac{1}{2}(xdw-wdx+zdy-ydz)\right)^2\\& \quad
+\left(ds+\frac{1}{2}(ydw-wdy+xdz-zdx)\right)^2
+\left(dt+\frac{1}{2}(zdw-wdz+ydx-xdy)\right)^2.
\end{split}
\end{equation*}

Then, by extending this metric, it seems natural to find a 7-dimensional generalization of the 3-dimensional $BCV$ spaces endowed with the two-parameter family of metrics

\begin{align*}
{\rm {ds}}^2_{m,l} & = \frac{dw^2+dx^2+dy^2+dz^2}{K^2} + \left(dr + \frac{l}{2} \frac{wdx-xdw+ydz-zdy}{K}\right)^2 \\
& \quad + \left(ds + \frac{l}{2} \frac{wdy-ydw+zdx-xdz}{K}\right)^2 +  \left(dt + \frac{l}{2} \frac{wdz-zdw+xdy-ydx}{K}\right)^2,
\end{align*}
where $l, m$ are real numbers and $K =1+m(w^2 +x^2+y^2+z^2)>0$.

Then $(EBCV, {\rm {ds}}^2_{m,l})$ will be called extended $BCV$ spaces ($EBCV$ for short).

That metric is obtained as a conformal deformation of the Euclidean metric
of $\mathbb{R}^4$ by adding three suitable terms which depend on $l$ and
$m$ concerning the imaginary part of $q \bar q$, for a quaternion $q$.
When $m=0$ we get a one-parameter of Riemannian metrics depending on $l$.
Furthermore, if $l=1$, we find the 7-dimensional quaternionic Heisenberg
group (see \cite{Kunkel} and \cite{W}). The manifold $EBCV$ provides
another example of sub-riemannian geometry and the horizontal distribution
is a 2-step breaking-generating distribution everywhere.

Observe that when $m=l=0$, $EBCV$ is nothing but $\mathbb{R}^7$; when $m>0, l=0$, $EBCV \cong \mathbb{S}^4(4m) \times \mathbb{R}^3$ and when $m<0, l=0$, $EBCV \cong \mathbb{H}^4(4m) \times \mathbb{R}^3$.


The metric ${\rm ds}^2_{m,l}$ can also be written as
$${\rm ds}^2_{m,l} =  \sum_{{\alpha}=1}^{7} \omega ^{\alpha} \otimes  \omega ^{\alpha}, $$
where

\begin{center}
\begin{tabular}{ll}

$\omega^1= dr + \frac{l}{2K}(wdx-xdw+ydz-zdy)$, \hspace{2cm} & $\omega^4=\frac{1}{K} dw$,\\[3mm]

$\omega^2= ds + \frac{l}{2K}(wdy-ydw+zdx-xdz)$, & $\omega ^5=\frac{1}{K} dx$,   \\[3mm]

$\omega^3= dt + \frac{l}{2K}(wdz-zdw+xdy-ydx)$, & $\omega ^6=\frac{1}{K} dy$,  \\[3mm]

  & $\omega^7=\frac{1}{K} dz$,

\end{tabular}
\end{center}

with the corresponding dual orthonormal frame

        $$X_1=\partial_r, \hspace{3cm} X_2=\partial_s, \hspace{3cm} X_3=\partial_t,$$

        $$X_4=K \partial_w+\frac {lx}{2} \partial_r +\frac {ly}{2} \partial_s + \frac {lz}{2}  \partial_t,\qquad X_5=K \partial_x-\frac {lw}{2} \partial_r -\frac {lz}{2} \partial_s + \frac {ly}{2}  \partial_t,$$

        $$X_6=K \partial_y+\frac {lz}{2} \partial_r -\frac {lw}{2} \partial_s - \frac {lx}{2}  \partial_t,\qquad
        X_7=K \partial_z-\frac {ly}{2} \partial_r +\frac {lx}{2} \partial_s - \frac {lw}{2}  \partial_t.$$

%

Writing $1\leq i, j \leq 3$, $4\leq a \leq 7$, we find that
$$ [X_i, X_j] = 0;\qquad   [X_i, X_a] = 0,$$
as well as
$$[X_4, X_5] = -l\{1+m(x^2+y^2)\}X_1+ml(wz+xy)X_2-ml(wy-xz)X_3-2mxX_4+2mwX_5,$$ and so on (see Appendix).


For later use, when $m=0$ brackets reduce to

\begin{center}
\begin{tabular}{lll}
$[X_4, X_5] = -l X_1$, \quad &  \quad $[X_4, X_6] =-l X_2$,  \quad & \quad  $[X_4, X_7] =-l X_3$,\\
$[X_5, X_6] =-l X_3$,  \quad & \quad  $[X_5, X_7] =l X_2$,  \quad & \quad  $[X_6, X_7] =-l X_1$.\\
\end{tabular}
\end{center}



\begin{remark}
When $l=1$, we have the brackets of the quaternionic contact manifold.
\end{remark}

As for the Levi-Civita connection we find out
$$\nabla_{X_i}X_j= 0,\qquad \nabla_{X_i}X_a= \nabla_{X_a}X_i,$$
and
\begin{equation*}\footnotesize{
\begin{split}
\nabla_{X_1}X_4 & = \frac{l}{2}\{1+m(y^2+z^2\}X_5+\frac{ml}{2}(wz-xy)X_6-\frac{ml}{2}(wy+xz)X_7,\\
\nabla_{X_1}X_5 & =-\frac{l}{2}\{1+m(y^2+z^2)\}X_4+\frac{ml}{2}(wy+xz)X_6+\frac{ml}{2}(wz-xy)X_7,\\
\nabla_{X_1}X_6 & =-\frac{ml}{2}(wz-xy)X_4-\frac{ml}{2}(wy+xz)X_5+\frac{l}{2}\{1+m(w^2+x^2)\}X_7,\\
\nabla_{X_1}X_7 & =\frac{ml}{2}(wy+xz)X_4-\frac{ml}{2}(wz-xy)X_5-\frac{l}{2}\{1+m(w^2+x^2)\}X_6,\\
\end{split}}
\end{equation*}
and son on (see Appendix).



\medskip

When $m=0$, the Levi-Civita connection reduces to

\begin{center}
\begin{tabular}{llll}
$\nabla_{X_1}X_4= \frac{l}{2} X_5$, \quad \quad & $\nabla_{X_3}X_4= \frac{l}{2} X_7$, \quad \quad &  $\nabla_{X_5}X_4= \frac{l}{2} X_1$, \quad \quad & $\nabla_{X_7}X_4=\frac{l}{2} X_3$,\\
$\nabla_{X_1}X_5=-\frac{l}{2} X_4$, \qquad & $\nabla_{X_3}X_5=\frac{l}{2} X_6$, \qquad &   $\nabla_{X_5}X_5=0$, \quad \quad & $\nabla_{X_7}X_5=-\frac{l}{2} X_2$,\\
$\nabla_{X_1}X_6=\frac{l}{2} X_7$, \quad \quad & $\nabla_{X_3}X_6-\frac{l}{2} X_5$, \quad \quad &   $\nabla_{X_5}X_6=-\frac{l}{2} X_3$, \quad \quad & $\nabla_{X_7}X_6=\frac{l}{2} X_1$,\\
$\nabla_{X_1}X_7=-\frac{l}{2} X_6$, \quad \quad & $\nabla_{X_3}X_7=-\frac{l}{2} X_4$, \quad \quad &   $\nabla_{X_5}X_7=\frac{l}{2} X_2$, \quad \quad & $\nabla_{X_7}X_7=0$.\\
$\nabla_{X_2}X_4= \frac{l}{2} X_6$, \quad \quad & $\nabla_{X_4}X_4= 0$, \quad \quad &   $\nabla_{X_6}X_4= \frac{l}{2} X_2$, \quad \quad & \\
$\nabla_{X_2}X_5=-\frac{l}{2} X_7$, \quad \quad & $\nabla_{X_4}X_5=-\frac{l}{2} X_1$, \quad \quad &   $\nabla_{X_6}X_5=\frac{l}{2} X_3$, \quad \quad & \\
$\nabla_{X_2}X_6=-\frac{l}{2} X_4$, \quad \quad & $\nabla_{X_4}X_6=-\frac{l}{2} X_2$, \quad \quad &   $\nabla_{X_6}X_6=0$, \quad \quad & \\
$\nabla_{X_2}X_7=\frac{l}{2} X_5$, \quad \quad & $\nabla_{X_4}X_7=-\frac{l}{2} X_3$, \quad \quad &   $\nabla_{X_6}X_7=-\frac{l}{2} X_1$, & \\
\end{tabular}
\end{center}


\begin{remark}
When $l=1$, we find the Levi-Civita connection of the quaternionic contact manifold.
\end{remark}

As for the curvature tensor $R$ we have

\begin{equation*}
\begin{split}
R_{X_1X_4X_1X_4}=R_{X_1X_5X_1X_5} & =\frac{l^2 }{4}\{1+m(K+1)(y^2+z^2)\},\\
R_{X_1X_6X_1X_6}=R_{X_1X_7X_1X_7} & =\frac{l^2 }{4}\{1+m(K+1)(w^2+x^2)\},\\
\end{split}
\end{equation*}
and so on (see Appendix).


\begin{remark}
When $m=0$, the curvature of the quaternionic contact manifold reduces to
\end{remark}
$$
R_{X_1X_4X_1X_4} = \frac{l^2}{4}, \quad  \ldots \quad R_{X_6X_7X_6X_7} =
-\frac{3l^2}{4}.
$$

\subsection{The Ricci tensor}

\begin{proposition}
	The matrix representing the Ricci tensor is given by
	
	\begin{center}
		{\footnotesize {\allowdisplaybreaks
				\begin{align*}
				&
				\left(\begin{array}{llll}
				\frac{l^2}{2}(K^2 +1) & 0 & 0 &\vdots \cr
				0 & \frac{l^2}{2}(K^2 +1)  & 0 & \vdots  \cr
				0 &0 & \frac{l^2}{2}(K^2 +1)  & \vdots  \cr
				-mlx(K+2) &-mly(K+2) &-mlz(K+2) & \vdots  \cr
				mlw(K+2) &mlz(K+2) &-mly(K+2) & \vdots  \cr
				-mlz(K+2) &mlw(K+2) &mlx(K+2) & \vdots  \cr
				mly(K+2) &-mlx(K+2) & mlw(K+2) & \vdots  \cr
				\end{array}
				\right.&\\
				& \left.\begin{array}{rrrrr}
				\vdots& -mlx(K+2) & mlw(K+2) &- mlz(K+2) &mly(K+2)  \cr
				\vdots& -mly(K+2) & mlz(K+2) & mlw(K+2) &- mlx(K+2)  \cr
				\vdots& -mlz(K+2) & -mly(K+2) & mlx(K+2) & mlw(K+2)  \cr
				\vdots& A(K-1-mw^2)+B& ml^2(K+1)wx& ml^2(K+1)wy & ml^2(K+1)wz  \cr
				\vdots& ml^2(K+1)wx & A(K-1-mx^2)+B& ml^2(K+1)xy & ml^2(K+1)xz  \cr
				\vdots& ml^2(K+1)wy &ml^2(K+1)xy & A(K-1-my^2)+B& ml^2(K+1)yz  \cr
				\vdots& ml^2(K+1)wz &ml^2(K+1)xz & ml^2(K+1)yz & A(K-1-mz^2)+B \cr
				\end{array}\right).&
				\end{align*}
		}}
	\end{center}
	where $A =-l^2(K+1)$ and $B = 12m-3/2l^2$.
\end{proposition}


Some particular cases could be interesting, for instance we get the following Ricci matrix when
$K=1$ (or $m=0$)

%

$${\rm Ric}_1=\left(  \begin{array} {ccccccc}
l^2& 0 & 0 &0&0&0&0 \cr
0& l^2 & 0 &0&0&0&0 \cr
0& 0 & l^2  &0&0&0&0  \cr
0& 0 & 0 &-3/2l^2 &0&0&0  \cr
0& 0 & 0 &0&-3/2l^2 &0&0  \cr
0& 0 & 0 &0&0&-3/2l^2 &0 \cr
0& 0 & 0 &0&0&0&-3/2l^2  \cr
\end{array}  \right)$$

%
%

%

\begin{remark}
	When $l=1$, we find the Ricci curvature of the quaternionic contact manifold.
\end{remark}

An easy computation leads to
\begin{corollary}
	The $EBCV$ manifold has constant scalar  curvature $S= 48 m$.
\end{corollary}


\section{The homogeneous structure}

In \cite{AS} W. Ambrose and I. M. Singer proved that a connected, complete and simply-connected Riemannian manifold $(M,g)$ is homogeneous if and only if there exists a (1,2) tensor field $T$ such that

(i) $g(T_XY,Z)+g(Y,T_XZ)=0$,

(ii) $(\nabla_XR)_{YZ}=[T_x,R_{YZ}]-R_{T_XYZ}-R_{YT_XZ}$,

(iii) $(\nabla_XT)_{Y}=[T_X,T_{Y}]-T_{T_XY}$,

\noindent for $X,Y,Z \in \mathfrak{X}(M)$, where $\nabla$ stands for the Levi-Civita connection and $R$ is the Riemann curvature tensor of $M$ (see \cite{TV}). As a consequence, Tricerri and Vanhecke define a homogeneous Riemannian structure on $(M,g)$ as a (1,2) tensor field $T$ which is a solution of the above three equations. Instead of taking (1,2) tensors it is prefered to work with (0,3) tensors via the isomorphism $T_{uvw}=g(T_uv,w)$, for $u, v, w \in T_pM$ and $p\in M$.

Then they consider the vector space $\mathfrak{T}$ of (0,3) tensors having the same symmetries as a homogeneous structure, i. e., $\mathfrak{T}=\{T : T_{uvw}=-T_{uwv},\; u, v, w \in T_pM\}$. The natural action of the orthogonal group $O(T_pM)$ on $\mathfrak{T}$ gives us the decomposition into eight irreducible invariant components. The main building blocks are defined as follows:

\begin{equation*}
\begin{split}
\mathfrak{T}_1&=\{T\in \mathfrak{T} : T_{uvw}=g(u,v)\alpha(w)-g(u,w)\alpha(v), \alpha \in T^*_p(M)\},\\
\mathfrak{T}_2&=\{T\in \mathfrak{T} : T_{uvw}+T_{wuv}+T_{vwu}=0, c_{12}(T)=0\},\\
\mathfrak{T}_3&=\{T\in \mathfrak{T} : T_{uvw}+T_{vuw}=0\},
\end{split}
\end{equation*}
where $u, v, w \in T_pM$ and $c_{12}(T)(w)=\sum_i T_{e_ie_iw}$, for any orthonormal basis $\{e_i\}$ of $T_pM$.

We consider in $EBCV$ the characteristic connection $D$ defined by (see \cite{Ro}):
$$D_{A}B =  \nabla_{A}B + \frac{P}{2} (\nabla_{A}P)B,$$
where $P$ is the natural almost product structure given by $P = \mathcal{V}-\mathcal{H}$, $Id=\mathcal{V}+\mathcal{H}$. Let us remember that the vertical distribution in $EBCV$ is spanned by $X_1, X_2, X_3$ and the horizontal distribution by $X_4, X_5, X_6, X_7$. Then we have
\begin{equation*}
\begin{split}
D_{X_i}X_j & = \mathcal{V} (\nabla_{X_i}X_j),\; i, j = 1, 2, 3,\\
D_{X_a}X_j & = \mathcal{V} (\nabla_{X_a}X_j),\; a= 4\ldots, 7;\;  j = 1, 2,3,\\
D_{X_i}X_b & = \mathcal{H} (\nabla_{X_i}X_b),\; i = 1, 2, 3;\; b= 4, \ldots,7,\\
D_{X_a}X_b & = \mathcal{H} (\nabla_{X_a}X_b),\; a, b = 4, \ldots,7.
\end{split}
\end{equation*}

This is a metric connection which makes parallel both the curvature and the torsion tensors and can be completely obtained by using the table giving the Levi-Civita connection.

By denoting $T^D$ the torsion tensor of $D$, that is,
$$T^D_L M\equiv T^D (L, M) = D_{L}M - D_{M}L - [L, M],$$ or equivalently
$$ T^D (L , M) = \frac{P}{2}\left((\nabla_LP)M - (\nabla_MP)L\right),$$
we find out
$$T^D (X_i, X_j) = T^D (X_i, X_a) = 0,\; i, j = 1, 2, 3;\; a = 4, \ldots, 7,$$

$$T^D (X_4, X_4) = T^D (X_5, X_5) = T^D (X_6, X_6)= T^D (X_7, X_7) = 0,$$
as well as

\begin{equation*}
\begin{split}
T^D (X_4, X_5) &= l\{(1+m(y^2 + z^2)) X_1 -m(wz +xy) X_2 + m(wy-xz)X_3\},\\
T^D (X_4, X_6) &= l\{m(wz-xy) X_1 + (1+m(x^2 + z^2)) X_2 - m(wx + yz)X_3\},\\
T^D (X_4, X_7) &= l \{-m(wy +xz) X_1 + (m(wx - yz) X_2 + (1+m(x^2 + y^2)X_3\},\\
T^D (X_5, X_6) &= l \{m(wy +xz) X_1 - (m(wx - yz) X_2 + (1+m(w^2 + z^2)X_3\},\\
T^D (X_5, X_7) &= l \{m(wz - xy) X_1 - (1+m(w^2 + y^2) X_2 - m(wx + yz)X_3\},\\
T^D (X_6, X_7) &= l \{(1+m(w^2 + x^2)) X_1 + m(wz + xy) X_2 - m(wy -xz)X_3\}.
\end{split}
\end{equation*}

Then $T^D$ defines a homogeneous structure on $EBCV$ in the sense of Tricerri-Vanhecke (see \cite{TV}, pags I and 15-16). Furthermore, it is easy to see that
$$c_{12} (T) = \sum T_{X_r} X_r = 0,\; r =1, \cdots, 7,$$
so that $T^D$ defines a homogeneous structure which is lying in the class $\mathfrak{T}_2 \oplus \mathfrak{T}_3$.

However, $T^D$ does not belong to  $\mathfrak{T}_2$, since, for instance,
$$T^D_{X_1 X_4 X_5} + T^D_{X_5 X_1 X_4} + T^D_{X_4 X_5 X_1} = \langle T^D_{X_1} X_4, X_5 \rangle + \langle T^D_{X_5} X_1, X_4 \rangle
+\langle T^D_{X_4} X_5, X_1 \rangle = l\{1+m(y^2 + z^2)\} \neq 0.$$
Finally, from the definition of $T^D$ we have that
$$T^D_{XYZ} = \langle T^D_{X} Y, Z \rangle = -\langle T^D_{Y} X, Z \rangle, $$ that is,
$T^D_{XYZ} + T^D_{YXZ} = 0$, and therefore $T^D$ is lying in $\mathfrak{T}_3$.

\newpage

\section{Killing vector fields in $EBCV$}

Remember that a Killing vector field  is a vector field on a Riemannian manifold that preserves the metric. Killing vector fields are the infinitesimal generators of isometries, that is, flows generated by Killing fields are continuous isometries of the manifold.
Specifically, a vector field $X$ is a Killing vector field if the Lie derivative with respect to $X$ of the metric g vanishes: $ {\displaystyle {\mathcal {L}}_{X}g=0\,} $ or equivalently

\begin{equation} \label{KVF1}
\mathcal {L}_X{\rm ds}^2_{l,m} = (\mathcal {L}_X \omega ^{\alpha})\otimes \omega^{\alpha} = 0,
\end{equation} where
$$\mathcal {L}_X \omega ^{\alpha} = \iota_Xd \omega ^{\alpha} + d(\iota_X \omega ^{\alpha}).$$

In terms of the Levi-Civita connection, Killing's condition is equivalent to

\begin{equation} \label{KVF2}
g(\nabla _{Y}X,Z)+g(Y,\nabla _{Z}X)=0.
\end{equation}

It is easy to prove that
\begin{proposition}   $\mathcal {L}_{X}g(Y,Z) =0$ if and only if
$\mathcal {L}_{X}g(X_i,X_j) =0$ for basic vector fields $X_i, X_j$.
\end{proposition}

We know that the dimension of the Lie algebra of the Killing vector fields is $m \leq n(n+1)/2$
and the maximum is reached on constant curvature manifolds (\cite{K-N}, p. 238, Vol. II) , then
for our manifold  $m < 28$. Then obviously

\begin{proposition}   The basic vertical vector fields $X_1, X_2, X_3$ are Killing fields.
\end{proposition}

From (\ref{KVF2}) it is easy to prove that
the horizontal basic vector fields $X_4,\cdots, X_7$ are not Killing vector fields.

 In her thesis, Profir \cite{Profir} proved that the Lie algebra of Killing vector fields of $BCV$ spaces is $4$-dimensional. Our problem now is to determine the space of Killing vector fields in $EBCV$.

\subsection{The Killing equations}

In the usual coordinate system $(r,s,t,w,x,y,z)$ on $EBCV$, a vector field $X=\sum_{{\alpha}=1}^{7}f_{\alpha}X_{\alpha}$
will be a Killing field if and only if the real functions $f_i$ satisfy the following system of $28$-partial differential equations:

%
%
%
%
%
%
%
%
%
%

\begin{center}
\def\arraystretch{1.3}
\begin{tabular}{l}
$\partial_r(f_1) = 0$,\\

$\partial_s(f_2) = 0$,\\

$\partial_t(f_3) = 0$,\\

$\partial_r(f_2)+\partial_s(f_1) = 0$,\\

$\partial_r(f_3)+\partial_t(f_1) = 0$,\\

$\partial_s(f_3)+\partial_t(f_2) = 0$,\\

\footnotesize{$\partial_r(f_4)+K\partial_w(f_1)+\frac{ly}{2}\partial_s(f_1)+\frac{lz}{2}\partial_t(f_1)-
l\{1+m(y^2 + z^2)\}f_5-ml(wz-xy)f_6+ml(wy+xz)f_7= 0$},\\

\footnotesize{$\partial_r(f_5)+K\partial_x(f_1)-\frac{lz}{2}\partial_s(f_1)+\frac{ly}{2}\partial_t(f_1)+
l\{1+m(y^2 + z^2)\}f_4-ml(wy+xz)f_6-ml(wz-xy)f_7= 0$},\\

\footnotesize{$\partial_r(f_6)+K\partial_y(f_1)-\frac{lw}{2}\partial_s(f_1)-\frac{lx}{2}\partial_t(f_1)+ml(wz-xy)f_4+
ml(wy+xz)f_5-l\{1+m(w^2+x^2\}f_7= 0$},\\

\footnotesize{$\partial_r(f_7)+K\partial_z(f_1)+\frac{lx}{2}\partial_s(f_1)-\frac{lw}{2}\partial_t(f_1)-ml(wy+xz)f_4+
ml(wz-xy)f_5+l\{1+m(w^2+x^2\}f_6= 0$},\\

\footnotesize{$\partial_s(f_4)+K\partial_w(f_2)+\frac{lx}{2}\partial_r(f_2)+\frac{lz}{2}\partial_t(f_2)+
ml(wz+xy)f_5-l\{1+m(x^2+z^2)\}f_6-ml(wx-yz)f_7= 0$},\\

\footnotesize{$\partial_s(f_5)+K\partial_x(f_2)-\frac{lw}{2}\partial_r(f_2)+\frac{ly}{2}\partial_t(f_2)-ml(wz+xy)f_4+
ml(wx-yz)f_6+l\{1+m(w^2+y^2)\}f_7= 0$},\\

\footnotesize{$\partial_s(f_6)+K\partial_y(f_2)-\frac{lz}{2}\partial_r(f_2)-\frac{lx}{2}\partial_t(f_2)+l\{1+m
(x^2+z^2)\}f_4-ml(wx-yz)f_5-ml(wz+xy)f_7= 0$},\\

\footnotesize{$\partial_s(f_7)+K\partial_z(f_2)-\frac{ly}{2}\partial_r(f_2)-\frac{lw}{2}\partial_t(f_2)+ml(wx-yz)f_4-
l\{1+m(w^2+y^2)\}f_5+ml(wz+xy)f_6= 0$},\\

\footnotesize{$\partial_t(f_4)+K\partial_w(f_3)+\frac{lx}{2}\partial_r(f_3)+\frac{ly}{2}\partial_s(f_3)+ml(wy-xz)f_5+
ml(wx+yz)f_6-l\{1+m(x^2 + y^2)\}f_7= 0$},\\

\footnotesize{$\partial_t(f_5)+K\partial_x(f_3)-\frac{lw}{2}\partial_r(f_3)-\frac{lz}{2}\partial_s(f_3)+ml(wy-xz)\}f_4-
l\{1+m(w^2+z^2)\}f_6+ml(wx+yz)f_7= 0$},\\

\footnotesize{$\partial_t(f_6)+K\partial_y(f_3)+\frac{lz}{2}\partial_r(f_3)-\frac{lw}{2}\partial_s(f_3)-ml(wx+yz)f_4+
l\{1+m(w^2+z^2\}f_5+ml(wy-xz)f_7= 0$},\\

\footnotesize{$\partial_t(f_7)+K\partial_z(f_3)-\frac{ly}{2}\partial_r(f_3)+\frac{lx}{2}\partial_s(f_3)+
l\{1+m(x^2+y^2)\}f_4-ml(wx+yz)f_5-ml(wy-xz)f_6= 0$},\\

\footnotesize{$K\partial_w(f_4)+\frac{lx}{2}\partial_r(f_4)+\frac{ly}{2}\partial_s(f_4)+\frac{lz}{2}\partial_t(f_4)-2mx f_5-2my f_6-2mz f_7=0$},\\

\footnotesize{$K\partial_w(f_5)+\frac{lx}{2}\partial_r(f_5)+\frac{ly}{2}\partial_s(f_5)+\frac{lz}{2}\partial_t(f_5)+
K\partial_x(f_4)-\frac{lw}{2}\partial_r(f_4)-\frac{lz}{2}\partial_s(f_4)+\frac{ly}{2}\partial_t(f_4)+2mx f_4+2mw f_5=0$},\\

\footnotesize{$K\partial_w(f_6)+\frac{lx}{2}\partial_r(f_6)+\frac{ly}{2}\partial_s(f_6)+\frac{lz}{2}\partial_t(f_6)+
K\partial_y(f_4)+\frac{lz}{2}\partial_r(f_4)-\frac{lw}{2}\partial_s(f_4)-\frac{lx}{2}\partial_t(f_4)+2my f_4+2mw f_6=0$},\\

\footnotesize{$K\partial_w(f_7)+\frac{lx}{2}\partial_r(f_7)+\frac{ly}{2}\partial_s(f_7)+\frac{lz}{2}\partial_t(f_7)+
K\partial_z(f_4)-\frac{ly}{2}\partial_r(f_4)+\frac{lx}{2}\partial_s(f_4)-\frac{lw}{2}\partial_t(f_4)+2mz f_4+2mw f_7=0$},\\

\footnotesize{$K\partial_x(f_5)-\frac{lw}{2}\partial_r(f_5)-\frac{lz}{2}\partial_s(f_5)+\frac{ly}{2}\partial_t(f_5)-2mw f_4-2my f_6 - 2mz f_7=0$},\\

\footnotesize{$K\partial_x(f_6)-\frac{lw}{2}\partial_r(f_6)-\frac{lz}{2}\partial_s(f_6)+\frac{ly}{2}\partial_t(f_6)+
K\partial_y(f_5)+\frac{lz}{2}\partial_r(f_5)-\frac{lw}{2}\partial_s(f_5)-\frac{lx}{2}\partial_t(f_5)+2my f_5+2mx f_6=0$},\\

\footnotesize{$K\partial_x(f_7)-\frac{lw}{2}\partial_r(f_7)-\frac{lz}{2}\partial_s(f_7)+\frac{ly}{2}\partial_t(f_7)+
K\partial_z(f_5)-\frac{ly}{2}\partial_r(f_5)+\frac{lx}{2}\partial_s(f_5)-\frac{lw}{2}\partial_t(f_5)+2mz f_5+2mx f_7=0$},\\

\footnotesize{$K\partial_y(f_6)+\frac{lz}{2}\partial_r(f_6)-\frac{lw}{2}\partial_s(f_6)-\frac{lx}{2}\partial_t(f_6)-2mw f_4-2mx f_5-2mz f_7=0$},\\

\footnotesize{$K\partial_y(f_7)+\frac{lz}{2}\partial_r(f_7)-\frac{lw}{2}\partial_s(f_7)-\frac{lx}{2}\partial_t(f_7)+
K\partial_z(f_6)-\frac{ly}{2}\partial_r(f_6)+\frac{lx}{2}\partial_s(f_6)-\frac{lw}{2}\partial_t(f_6)+2mz f_6+2my f_7=0$},\\

\footnotesize{$K\partial_z(f_7)-\frac{ly}{2}\partial_r(f_7)+\frac{lx}{2}\partial_s(f_7)-\frac{lw}{2}\partial_t(f_7)-2mw f_4-2mx f_5-2my f_6=0$}.\\
\end{tabular}
\end{center}

\newpage

It seems that the solution of the system is very difficult, so that we focus on solving the system for $m=0$, that is:

%
%
%
%
%
%
%
%
%
%

\begin{center}
\def\arraystretch{1.3}
\begin{tabular}{l}

$\partial_r(f_1)=0$,\\

$\partial_s(f_2)=0$,\\

$\partial_t(f_3)=0$,\\

$\partial_r(f_2)+\partial_s(f_1)=0$,\\

$\partial_r(f_3)+\partial_t(f_1)=0$,\\

$\partial_s(f_3)+\partial_t(f_2)=0$,\\

$\partial_r(f_4)+\partial_w(f_1)+\frac{ly}{2}\partial_s(f_1)+\frac{lz}{2}\partial_t(f_1)-lf_5=0$,\\

$\partial_r(f_5)+\partial_x(f_1)-\frac{lz}{2}\partial_s(f_1)+\frac{ly}{2}\partial_t(f_1)+lf_4=0$,\\

$\partial_r(f_6)+\partial_y(f_1)-\frac{lw}{2}\partial_s(f_1)-\frac{lx}{2}\partial_t(f_1)-lf_7=0$,\\

$\partial_r(f_7)+\partial_z(f_1)+\frac{lx}{2}\partial_s(f_1)-\frac{lw}{2}\partial_t(f_1)+lf_6=0$,\\

$\partial_s(f_4)+\partial_w(f_2)+\frac{lx}{2}\partial_r(f_2)+\frac{lz}{2}\partial_t(f_2)-lf_6=0$,\\

$\partial_s(f_5)+\partial_x(f_2)-\frac{lw}{2}\partial_r(f_2)+\frac{ly}{2}\partial_t(f_2)+lf_7=0$,\\

$\partial_s(f_6)+\partial_y(f_2)+\frac{lz}{2}\partial_r(f_2)-\frac{lx}{2}\partial_t(f_2)+lf_4=0$,\\

$\partial_s(f_7)+\partial_z(f_2)-\frac{ly}{2}\partial_r(f_2)-\frac{lw}{2}\partial_t(f_2)-lf_5=0$,\\

$\partial_t(f_4)+\partial_w(f_3)+\frac{lx}{2}\partial_r(f_3)+\frac{ly}{2}\partial_s(f_3)-lf_7=0$,\\

$\partial_t(f_5)+\partial_x(f_3)-\frac{lw}{2}\partial_r(f_3)-\frac{lz}{2}\partial_s(f_3)-lf_6=0$,\\

$\partial_t(f_6)+\partial_y(f_3)+\frac{lz}{2}\partial_r(f_3)-\frac{lw}{2}\partial_s(f_3)+lf_5=0$,\\

$\partial_t(f_7)+\partial_z(f_3)-\frac{ly}{2}\partial_r(f_3)+\frac{lx}{2}\partial_s(f_3)+lf_4=0$,\\

$\partial_w(f_4)+\frac{lx}{2}\partial_r(f_4)+\frac{ly}{2}\partial_s(f_4)+\frac{lz}{2}\partial_t(f_4)=0$,\\

$\partial_w(f_5)+\frac{lx}{2}\partial_r(f_5)+\frac{ly}{2}\partial_s(f_5)+\frac{lz}{2}\partial_t(f_5)+\partial_x(f_4)-\frac{lw}{2}\partial_r(f_4)-\frac{lz}{2}\partial_s(f_4)+\frac{ly}{2}\partial_t(f_4)=0$,\\

$\partial_w(f_6)+\frac{lx}{2}\partial_r(f_6)+\frac{ly}{2}\partial_s(f_6)+\frac{lz}{2}\partial_t(f_6)+\partial_y(f_4)+\frac{lz}{2}\partial_r(f_4)-\frac{lw}{2}\partial_s(f_4)-\frac{lx}{2}\partial_t(f_4)=0$,\\

$\partial_w(f_7)+\frac{lx}{2}\partial_r(f_7)+\frac{ly}{2}\partial_s(f_7)+\frac{lz}{2}\partial_t(f_7)+\partial_z(f_4)-\frac{ly}{2}\partial_r(f_4)+\frac{lx}{2}\partial_s(f_4)-\frac{lw}{2}\partial_t(f_4)=0$,\\

$\partial_x(f_5)-\frac{lw}{2}\partial_r(f_5)-\frac{lz}{2}\partial_s(f_5)+\frac{ly}{2}\partial_t(f_5)=0$,\\

$\partial_x(f_6)-\frac{lw}{2}\partial_r(f_6)-\frac{lz}{2}\partial_s(f_6)+\frac{ly}{2}\partial_t(f_6)+\partial_y(f_5)+\frac{lz}{2}\partial_r(f_5)-\frac{lw}{2}\partial_s(f_5)-\frac{lx}{2}\partial_t(f_5)=0$,\\

$\partial_x(f_7)-\frac{lw}{2}\partial_r(f_7)-\frac{lz}{2}\partial_s(f_7)+\frac{ly}{2}\partial_t(f_7)+\partial_z(f_5)-\frac{ly}{2}\partial_r(f_5)+\frac{lx}{2}\partial_s(f_5)-\frac{lw}{2}\partial_t(f_5)=0$,\\

$\partial_y(f_6)+\frac{lz}{2}\partial_r(f_6)-\frac{lw}{2}\partial_s(f_6)-\frac{lx}{2}\partial_t(f_6)=0$,\\

$\partial_y(f_7)+\frac{lz}{2}\partial_r(f_7)-\frac{lw}{2}\partial_s(f_7)-\frac{lx}{2}\partial_t(f_7)+\partial_z(f_6)-\frac{ly}{2}\partial_r(f_6)+\frac{lx}{2}\partial_s(f_6)-\frac{lw}{2}\partial_t(f_6)=0$,\\

$\partial_z(f_7)-\frac{ly}{2}\partial_r(f_7)+\frac{lx}{2}\partial_s(f_7)-\frac{lw}{2}\partial_t(f_7)=0$.\\

\end{tabular}
\end{center}

\newpage

\noindent whose solution is given by

\begin{equation*}
\begin{split}
f_1(r,s,t,w,x,y,z) & =(P+R)s+(S-N)t+\frac{l}{2}\{-M(w^2+x^2)-U(y^2+z^2)+(R-P)(wy+xz)\\& \hspace{4cm}+(N+S)(wz-xy)+2Tw-2Qx+2Wy-2Vz\}+C_1,\\
f_2(r,s,t,w,x,y,z) & =-(P+R)r+(M+U)t-\frac{l}{2}\{N(w^2+y^2)-S(x^2+z^2)+(R-P)(wx-yz)\\& \hspace{4cm}+(M-U)(wz+xy)-2Vw+2Wx+2Qy-2Tz\}+C_2,\\
f_3(r,s,t,w,x,y,z) & =-(S-N)r-(M+U)s-\frac{l}{2}\{P(w^2+z^2)+R(x^2+y^2)+(N+S)(wx+yz)\\& \hspace{4cm}+(U-M)(wy-xz)-2Ww-2Vx+2Ty+2Qz\}+C_3,\\
f_4(r,s,t,w,x,y,z) & =Mx+Ny+Pz+Q,\\
f_5(r,s,t,w,x,y,z) & =-Mw+Ry+Sz+T,\\
f_6(r,s,t,w,x,y,z) & =-Nw-Rx+Uz+V,\\
f_7(r,s,t,w,x,y,z) & =-Pw-Sx-Uy+W,\\
\end{split}
\end{equation*}
where $M, N, P, Q, R, S, T, U, V, W, C_1, C_2, C_3 \in \mathbb{R}$.

\vspace{.5cm}

As a consequence, when $m=0$, we obtain

\begin{proposition}
The Lie algebra of Killing vector fields is 13-dimensional.	
\end{proposition}

%
%
%
%
%
%
%

%

\section{Computing horizontal geodesics of the quaternionic Hei-senberg group (see \cite{Mont})}

A Riemannian metric on a manifold $M$ is defined by a covariant two-tensor, which is to say, a section of the bundle $S^2(T^*M)$. There is no such object in subriemannian geometry. Instead, a subriemannian metric can be encoded as a contravariant symmetric two-tensor,  which is a section of $S^2(TM)$.  This two-tensor  has rank $k < n$, where $k$ is the rank of the distribution, so it cannot be inverted to obtain a Riemannian metric. We call this contravariant tensor the {\it cometric}.

\begin{definition} A {\rm cometric} is a section of the bundle $S^2(TM) \subset TM \otimes TM$ of symmetric bilinear  forms on the cotangent  bundle of $M$.
\end{definition}

Since $TM$ and $T^*M$ are dual,  any  cometric  defines  a fiber-bilinear  form
$((\cdot \,, \cdot)): T^*M \otimes T^*M \to \mathbb{R}$, i.e. a kind of inner product on covectors.  This form in turn defines a symmetric bundle  map $\beta: T^*M \to TM$ by $p(\beta_q(\mu)) = ((p, \mu))_q$,  for $p, \mu \in T_q^*M$ and $q \in M$. Thus $\beta_q(\mu) \in T_qM$.  The  adjective  symmetric   means  that $\beta$ equals  its adjoint $\beta^*: T^*M \to T^{**}M = TM$.

The cometric $\beta$ for a subriemannian geometry is uniquely defined by the following conditions:

\newpage

\begin{itemize}
\item[(1)] ${\rm im}(\beta_q) = \mathcal{H}_q$;
\item[(2)] $p( v) = \langle \beta_q(p),v \rangle$, for $v \in\mathcal{H}_q$, $p \in T_qM$,
\end{itemize}
where $\langle \beta_q(p),v \rangle_q$ is the subriemannian  inner  product  on $\mathcal{H}_q$. Conversely, any cometric of constant rank defines a subriemannian geometry whose underlying distribution has that rank.

\begin{definition} The fiber-quadratic  function $H(q, p)= \frac{1}{2}(p, p)_q$, where $(\cdot \,, \cdot)_q$ is the cometric on the fiber $T_q^*M$, is called the  {\rm subriemannian Hamiltonian},  or  the {\rm kinetic energy}.
\end{definition}

The Hamiltonian $H$ is related to length and energy as follows.  Suppose that $\gamma$
is a horizontal curve. Then, $\dot \gamma(t)=\beta_{\gamma(t)}(p)$, for same covector $p \in T_{\gamma(t)}^*M$, and $$\frac{1}{2}||\dot \gamma||^2 = H (q, p).$$
$H$ uniquely determines $\beta$ by polarization, and $\beta$ uniquely determines the subriemannian  structure. This proves the  following  proposition:

\begin{proposition} The subriemannian structure is uniquely determined by its Hamiltonian. Conversely, any nonnegative fiber-quadratic Hamiltonian of constant fiber rank $k$ gives  rise to a sub-riemannian structure whose underlying distribution has rank $k$ .
\end{proposition}

To compute the subriemannian  Hamiltonian we can start with a local frame $\{X_a\}, a = 1, \ldots , k$, of  vector  fields for $\mathcal{H}$. Thinking of the $X_a$ as fiber-linear functions on the cotangent bundle, we rename them $P_a$ so that
$$P_a(q, p)=p(X_a(q)),\quad q \in M, p \in T_q^*M.$$

\begin{definition} Let $X$ be a vector field on the manifold $M$. The fiber-linear function on the cotangent bundle $P_X:T^*M \to \mathbb{R}$, defined by $P_X(q, p) = p(X( q))$ is called the {\rm momentum function} for $X$.
\end{definition}

Thus the $P_a = P_{X_a}$ are the momentum functions for our horizontal frame.
If $X_a = \sum X_a^{\alpha}(x)\frac{\partial}{\partial x^{\alpha}}$ is the expression for $X_a$  relative to coordinates $x^{\alpha}$, then $P_{X_a}(x,p )=\sum X_a^{\alpha}(x)p_{\alpha}$, where $p_{\alpha} =P_{\frac{\partial}{\partial x^{\alpha}}}$ are the momentum functions for the coordinate vector fields. The $x^{\alpha}$ and $p_{\alpha}$ together form a coordinate system on $T^*M$. They are called {\it canonical coordinates}.

Let $g_{ab}(q) = \langle X_a(q), X_b(q)\rangle_q$ be the matrix of inner products defined by our horizontal frame. Let $g^{ab}(q)$ be its inverse matrix. Then $g^{a b}$ is a $k \times k$ matrix-valued function defined in some open set of $M$.

\begin{proposition} Let $P_a$ and $g^{ab}$ be the functions on $T^*M$ that are induced by a local horizontal frame ${X_a}$ as just  described. Then
\begin{equation}\label{hamiltonian1}
H(q,p)=\frac{1}{2}\sum g^{ab}(q) P_a(q,p)P_b(q,p).
\end{equation}
\end{proposition}

Indeed,
\begin{equation*}
\begin{split}
H(q,p)&=\frac{1}{2}(p,p)_q=\frac{1}{2}(\sum p_adx^a, \sum p_bdx^b)=
\frac{1}{2}\sum g^{ab}(q)(p_a, p_b)\\&=\frac{1}{2}\sum g^{ab}(q)(p(X_a)(q), p(X_b)(q))=
\frac{1}{2}\sum g^{ab}(q)P_a(q,p) P_b(q,p).
\end{split}
\end{equation*}

Note, in particular, that if the $X_a$ are an orthonormal frame for $\mathcal{H}$ relative to the subriemannian inner product, then $H = \frac{1}{2}P_a^2$.

{\bf Normal geodesics}. Like any smooth function on the cotangent bundle, our function $H$ generates a system of Hamiltonian differential equations. In terms  of canonical coordinates $(x^{\alpha},p_{\alpha})$, these differential equations are
\begin{equation}\label{normalgeodesicequations}
\dot x^{\alpha}=\dfrac{\partial H}{\partial p_{\alpha}}, \quad \dot p_{\alpha}=-\dfrac{\partial H}{\partial x^{\alpha}}.
\end{equation}

\begin{definition} The Hamiltonian differential equations {\rm (\ref{normalgeodesicequations})} are called the {\rm normal geodesic equations}.
\end{definition}

Riemannian geometry can be viewed as a special case of subriemannian geometry, one in which the distribution is the entire tangent bundle. The cometric is the usual inverse metric, written $g^{ij}$ in coordinates. The normal geodesic equations in the Riemannian case are simply the standard geodesic equations, rewritten on the cotangent bundle.

\subsection{One remark about the Hamiltonian}

We follow word for word the computations in \cite{Mont}. The vector fields
\begin{equation*}
\begin{split}
W & = \partial_w+\frac{1}{2}(x \partial_r+y \partial_s+z \partial_t),\\
X & = \partial_x-\frac{1}{2}(w \partial_r+z \partial_s-y \partial_t),\\
Y & = \partial_y+\frac{1}{2}(z \partial_r-w \partial_s-x \partial_t),\\
Z & = \partial_z-\frac{1}{2}(y \partial_r-x \partial_s-w \partial_t),\\
\end{split}
\end{equation*}
which are the old $X_4, \cdots, X_7$ ones, provided $m=0, l=1$, along with $\{\partial_r, \partial_s, \partial_t\}$, form an orthonormal frame for the quaternionic contac manifold $\mathbb{H} \times {\rm Im}\mathbb{H}$. This means that $\{W, X, Y Z\}$ frame the fourth plane $\mathcal{H}$ and they are orthonormal with respect to the inner product ${\rm ds}^2 = (dw^2 + dx^2 + dy^2 + dz^2)|_{\mathcal{H}}$ on the distribution.  According to the above discussion, the subriemannian Hamiltonian is
\begin{equation}\label{hamiltonianquaternionic}
H = \frac{1}{2}(P_W^2 + P_X^2 + P_Y^2 + P_Z^2),
\end{equation}
where $P_W, P_X,P_Y, P_Z$ are the momentum  functions of the vector fields $W, X, Y, Z$, respectively. Thus
\begin{equation*}
\begin{split}
P_W &= p_w + \frac{1}{2}(x p_r + y p_s + z p_t),\\
P_X &= p_x  - \frac{1}{2}(w p_r +z p_s -y p_t),\\
P_Y &= p_y + \frac{1}{2}(z p_r  - w p_s - x p_t),\\
P_Z &= p_z -  \frac{1}{2}(y p_r -x p_s +w p_t),\\
\end{split}
\end{equation*}
where $p_w, p_x, p_y, p_z, p_r, p_s, p_t$ are the fiber coordinates on the cotangent bundle of $\mathbb{R}^7$ corresponding to the cartesian coordinates $w, x, y, z, r, s, t$ on $\mathbb{R}^7$. Again, these fiber coordinates are defined by writing a covector as
$p = p_w dw + p_x dx + p_y dy + p_z dz + p_r dr + p_s ds + p_t dt$.
 Together, $(w, x, y, z, r, s, t, p_w, p_x, p_y, p_z, p_r, p_s, p_t)$ are global coordinates on the cotangent bundle $T^*\mathbb{R}^7 = \mathbb{R}^7 \oplus \mathbb{R}^7$. Hamilton's equations can be written \begin{equation}\label{hamiltonian4}
 \frac{df}{du} = \{f, H\},\quad f \in C^{\infty}(T^*\mathbb{R}^7),
 \end{equation}
which holds for any smooth function $f$. The function $H$ defines a vector field $X_H$, called the Hamiltonian vector field, which has a flow  $\Phi_u : T^*\mathbb{R}^7  \to  T^*\mathbb{R}^7$. Let $f : T^*\mathbb{R}^7 \to \mathbb{R}$ be any smooth function on the cotangent bundle. Form the $u$-dependent function $f_u = \Phi_u^* f$ by pulling $f$ back via the flow. Thus
$f_u(w, x, y, z, r, s, t, p_w, p_x, p_y, p_z, p_r, p_s, p_t) = f ( \Phi_u(w, x, y, z, r, s, t, p_w, p_x, p_y, p_z, p_r, p_s, p_t))$.
In other words
$\frac{df}{du} = X_H[f_u]$, which gives meaning to the left-hand side of Hamilton's equations.

    To define the right hand side, which is to say the vector field $X_H$, we will need the Poisson bracket. The Poisson bracket on the cotangent bundle $T^*\mathbb{R}^7$ of a manifold $\mathbb{R}^7$ is a canonical Lie algebra structure defined on the vector space $C^{\infty}(T^*\mathbb{R}^7)$ of smooth functions on $T^*\mathbb{R}^7$. The Poisson bracket is denoted $\{\cdot \,, \cdot\} : C^{\infty}  \times  C^{\infty} \to C^{\infty}$, where $C^{\infty}= C^{\infty}(T^*\mathbb{R}^7)$, and can be defined by the coordinate formula
    $$\{f, g\} =\frac{\partial f}{\partial x^i}\frac{\partial g}{\partial p_i}-\frac{\partial g}{\partial x^i}\frac{\partial f}{\partial p_i}.$$
This formula is valid in any canonical coordinate system, and can be shown to be coordinate independent. The Poisson bracket satisfies the Leibniz identity
    $$\{f, g h\} = g \{f, h\} + h \{f, g\},$$
which means that the operation $\{., H\}$ defines a vector field $X_H$, called the Hamiltonian vector field.
By letting the functions $f$ vary over the collection of coordinate functions $x^i$ and we get the more common form of Hamilton's equations
        $${\dot x}^{\alpha} =\frac{\partial H}{\partial p_{\alpha}}, \quad {\dot p}_{\alpha} =-\frac{\partial H}{\partial x^{\alpha}}.$$
Indeed, for the first one we take $f = w$ and $g = H$. Then
    $\{w, H\}=\frac{\partial w}{\partial x^i}\frac{\partial H}{\partial p_i}-\frac{\partial H}{\partial x^i}\frac{\partial w}{\partial p_i}$ if and only if $\dot w = \frac{\partial H}{\partial p_w}$.
Also we have
    $$\dot x = \frac{\partial H}{\partial p_x}, \quad \dot y = \frac{\partial H}{\partial p_y}, \quad \dot z = \frac{\partial H}{\partial p_z}.$$

These equations are in turn equivalent to the above formulation (\ref{hamiltonian4}), which is more convenient to use, because the momentum function $W \mapsto P_W$ is a Lie algebra anti-homomorphism from the Lie algebra of all smooth vector fields on $\mathbb{R}^7$ to $C(T^*\mathbb{R}^7)$ with the Poisson brackets:
\begin{equation}\label{poissonbracket2}
\begin{split}
\{P_W, P_X\} = &-P_{[W, X]},\quad \{P_W, P_Y\} = -P_{[W, Y]}, \quad \{P_W, P_Z\} = -P_{[W, Z]},\\
\{P_X, P_Y\} =& -P_{[X, Y]},\quad  \{P_X, P_Z\} = -P_{[X, Z]},\quad \{P_Y, P_Z\} = -P_{[Y, Z]}.
\end{split}
\end{equation}

Since all calculations are similar, we only prove the first one:

$$\{P_W, P_X\} = \{p_w + \frac{x}{2} p_r + \frac{y}{2} p_s + \frac{z}{2} p_t, p_x -\frac{w}{2} p_r -\frac{z}{2} p_s + \frac{y}{2} p_t\}= p_r = - P_{[W, X]}.$$

For the quaternionic contact group, with our choose of $W, X, Y, Z$ as a frame for $\mathcal{H}$, we compute

$$[W, X]=-{\partial_r},\quad [W, Y] =-{\partial_s},\quad [W, Z]=-{\partial_t},$$

$$[X, Y] =-\partial_t,\quad [X, Z] = \partial_s,\quad [Y, Z] =-\partial_t,$$

$$[W, \partial_r] = [W, \partial_rs] = [W, \partial_t] = [X, \partial_r] = [X, \partial_s] = [X, \partial_t]=0,$$

$$[Y, \partial_r] = [Y, \partial_s] = [Y, \partial_t] = [Z, \partial_r] = [Z, \partial_s] = [Z, \partial_t] = 0.$$

Thus
$$\{P_W, P_X\}= \partial_r:= P_r,\quad
\{P_W, P_Y\}= \partial_s:= P_s,\quad
\{P_W, P_Z\}= \partial_t:= P_t,$$

$$\{P_X, P_Y\}= P_t,\quad \{P_X, P_Z\}= - p_s = - P_s,\quad \{P_Y, P_Z\}= p_r = P_r$$

We can prove that
$$\{P_W, P_r\}= \{P_W, P_s\} = \{P_W, P_t\} = \{P_X, P_r\}= \{P_X, P_s\} = \{P_X, P_t\} = 0,$$
$$\{P_Y, P_r\}= \{P_Y, P_s\} = \{P_Y, P_t\} = \{P_Z, P_r\}= \{P_Z, P_s\} = \{P_Z, P_t\} = 0.$$

These relations can also easily be computed by hand, from our formulae for $P_W, P_X, P_Y, P_Z$ and the bracket in terms of $w, x, y, z, r, s, r, p_w, p_x, p_y, p_z, p_r, p_s, p_t$.

\begin{lemma} By letting $f$ vary over the functions $w, x, y, z, r, s, r, P_W, P_X, P_Y, P_Z, P_r, P_s, P_t$, using the bracket relations and equation {\rm (\ref{poissonbracket2})}, we find that Hamilton's equations are equivalent to the system

\begin{center}
\begin{tabular}{lllll}
$\dot w =P_W$, &  &  &  &  ${\dot P}_W  =  p_r P_X + p_s P_Y +  p_t P_Z$,\\
$\dot x =P_X$, &  &  &  &  ${\dot P}_X  =  -p_r P_W - p_s P_Z +  p_t P_Y$,\\
$\dot y =P_Y$, &  &  &  &  ${\dot P}_Y  = p_r P_Z -p_s P_W -p_t P_X$, \\
$\dot z =P_Z$, &  &  &  &  ${\dot P}_Z  =  - p_r P_Y + p_s P_X - p_t P_W$,\\
$\dot r  = \frac{1}{2} (x P_W -w P_X + z P_Y-y P_Z)$, &  &  &  &  ${\dot P}_r  = 0$,\\
$\dot s = \frac{1}{2} (y P_W -z P_X + x P_Y-w P_Z)$, &  &  &  &  ${\dot P}_s  = 0$,\\
$\dot t =  \frac{1}{2} (z P_W +y P_X - x P_Y-w P_Z)$, &  &  &  &  ${\dot P}_t  = 0$.\\
\end{tabular}
\end{center}
\end{lemma}

To see it, remember that $H= \frac{1}{2}(P_W^2 +P_X^2 + P_Y^2 + P_z^2 )$. Then
\begin{equation*}
\begin{split}
\dot w &= \{w, H\}= P_w \frac{\partial P_W}{\partial p_w} = P_W,\\
\dot x &= \{x, H\}= P_X \frac{\partial P_X}{\partial p_x}= P_X,\\
\dot y &= P_Y,\\
\dot z &= P_Z.
\end{split}
\end{equation*}

Also, considering that:

$$\frac{\partial P_W}{\partial p_r} =\frac{x}{2},\qquad \frac{\partial P_W}{\partial p_s} =\frac{y}{2},\qquad \frac{\partial P_W}{\partial p_t} =\frac{z}{2},$$

$$\frac{\partial P_X}{\partial p_r} =-\frac{w}{2},\qquad \frac{\partial P_X}{\partial p_s} =-\frac{z}{2},\qquad \frac{\partial P_X}{\partial p_t} =\frac{y}{2},$$

$$\frac{\partial P_Y}{\partial p_r} =\frac{z}{2},\qquad \frac{\partial P_Y}{\partial p_s} =-\frac{w}{2},\qquad \frac{\partial P_Y}{\partial p_t} =-\frac{x}{2},$$

$$\frac{\partial P_Z}{\partial p_r} =-\frac{y}{2},\qquad \frac{\partial P_Z}{\partial p_s} =\frac{x}{2},\qquad \frac{\partial P_Z}{\partial p_t} =-\frac{w}{2},$$

we have
$$\dot r =\frac{1}{2} (x P_W -w P_X + z P_Y-y P_Z).$$

Indeed,
\begin{equation*}
\begin{split}
\dot r &= \{r, H\} = P_W \frac{\partial P_W}{\partial p_r} + P_X \frac{\partial P_X}{\partial p_r} + P_Y \frac{\partial P_Y}{\partial p_r} + P_Z \frac{\partial P_Z}{\partial p_r}
\\&= \frac{1}{2} (x P_W -w P_X + z P_Y-y P_Z)\\
\dot s &= \{s, H\} = P_W \frac{\partial P_W}{\partial p_s} + P_X \frac{\partial P_X}{\partial p_s} + P_Y \frac{\partial P_Y}{\partial p_s} + P_Z \frac{\partial P_Z}{\partial p_s}
\\&= \frac{1}{2} (y P_W -z P_X + x P_Y-w P_Z)\\
\dot t &= \{t, H\} = P_W \frac{\partial P_W}{\partial p_t} + P_X \frac{\partial P_X}{\partial p_t} + P_Y \frac{\partial P_Y}{\partial p_t} + P_Z \frac{\partial P_Z}{\partial p_t}
\\&= \frac{1}{2} (z P_W +y P_X - x P_Y-w P_Z).\\
\end{split}
\end{equation*}

Working as above we obtain
\begin{equation*}
\begin{split}
{\dot P}_W &= \{P_W, H\} = p_r P_X + p_s P_Y +  p_t P_Z,\\
{\dot P}_X &= \{P_X, H\} = -p_r P_W - p_s P_Z +  p_t P_Y,\\
{\dot P}_Y &= \{P_Y, H\} = p_r P_Z -p_s P_W -p_t P_X,\\
{\dot P}_Z &= \{P_Z, H\} = - p_r P_Y + p_s P_X - p_t P_W.
\end{split}
\end{equation*}

Then we are ready to show the following

\begin{theorem} The horizontal geodesics of the quaternionic Heisenberg group are exactly the horizontal lifts of arcs of circles, including line segments as a degenerate case.
\end{theorem}

\noindent {\bf Proof.} It is not difficult to see that ${\dot P}_r = {\dot P}_s ={\dot P}_t = 0$. These equations assert that $P_r = p_r$, $P_s = p_s$ and $P_t = p_t$ are constant. The variables $r, s, t$ appears nowhere in the right-hand sides of these equations. It follows that the variables $w, x, y, z, P_W, P_X, P_Y, P_Z$ evolve independently of $r, s, t$, and so we can view the system as defining a one-parameter family of dynamical systems on $\mathbb{R}^8$ parameterized by the constant value of $P_r, P_s, P_t$.

Combine $w, x, y, z$ into a single quaternionic variable $ \omega = w + i x + j y + k z$ and taking into account the fourteen equations one has
$$\frac{d\omega}{du} = P_W + i P_X + j P_Y + k P_Z$$
The $u$-derivative of $P_W + i P_X + j P_Y + k P_Z$ is
$-(i p_r + j p_s + k p_t)(P_W + i P_X + j P_Y + k P_Z)$.
Then we have $\frac{d^2 \omega}{du^2} = -(i p_r + j p_s + k p_t)\frac{d\omega}{du}$, where $p_r$, $p_s$ and $p_t$ are constant.

By integrating the above expression we get
 $$P_W + i P_X + jP_Y +kP_Z= P(0){\rm exp}(-(i p_r + j p_s + k p_t) t ),$$ where  $P(0) = P_W(0) + i P_X(0) + j P_Y(0) + k P_Z(0)$ .

A second integration yields the general form of the geodesics on the quaternionic contact group:

\begin{equation*}
\begin{split}
\omega (u) &= w(u) + i x(u) + j y(u) + k z(u) =\\&
\frac{P(0)}{i p_r + j p_s + k p_t}\left({\rm exp}(-(i p_r + j p_s + k p_t)t-1) + w(0) + i x(0) + j y(0) + k z(0)\right),\\
r(u) &= r(0) +\frac{1}{2}\int_0^t {\rm Im}_I (\bar \omega \, d\omega),\\
s(u) &= s(0) +\frac{1}{2}\int_0^t {\rm Im}_J (\bar \omega \, d\omega),\\
r(u) &= t(0) +\frac{1}{2}\int_0^t {\rm Im}_K (\bar \omega \, d\omega).\\
\end{split}
\end{equation*}

\newpage

\section{Appendix}

The brackets:
\begin{equation*}
\begin{split}
[X_4, X_5] & = -l\{1+m(x^2+y^2)\}X_1+ml(wz+xy)X_2-ml(wy-xz)X_3-2mxX_4+2mwX_5,\\
[X_4, X_6] & =-ml(wz-xy)X_1-l\{1+m(x^2+z^2)\}X_2+ml(wx+yz)X_3-2myX_4+2mwX_6,\\
[X_4, X_7] & =ml(wy+xz)X_1-ml(wx-yz)X_2-l\{1+(x^2+y^2)\}X_3-2mzX_4+2mwX_7,\\
[X_5, X_6] & =-ml(wy+xz)X_1+ml(wx-yz)X_2-l\{1+m(w^2+z^2)\}X_3-2myX_5+2mxX_6,\\
[X_5, X_7] & =ml(xy-wz)X_1+ l\{1+m(w^2+y^2)\}X_2+ml(wx+yz)X_3-2mzX_5+2mxX_7,\\
[X_6, X_7] & =-l\{1+m(w^2+x^2)\}X_1-ml(wz+xy)X_2+ml(wy-xz)X_3-2mzX_6+2myX_7.\\
\end{split}
\end{equation*}
The Levi-Civita connection:
\begin{equation*}\footnotesize{
\begin{split}
\nabla_{X_1}X_4 & = \frac{l}{2}\{1+m(y^2+z^2)\}X_5+\frac{ml}{2}(wz-xy)X_6-\frac{ml}{2}(wy+xz)X_7,\\
\nabla_{X_1}X_5 & =-\frac{l}{2}\{1+m(y^2+z^2)\}X_4+\frac{ml}{2}(wy+xz)X_6+\frac{ml}{2}(wz-xy)X_7,\\
\nabla_{X_1}X_6 & =-\frac{ml}{2}(wz-xy)X_4-\frac{ml}{2}(wy+xz)X_5+\frac{l}{2}\{1+m(w^2+x^2)\}X_7,\\
\nabla_{X_1}X_7 & =\frac{ml}{2}(wy+xz)X_4-\frac{ml}{2}(wz-xy)X_5-\frac{l}{2}\{1+m(w^2+x^2)\}X_6,\\
\nabla_{X_2}X_4 & = -\frac{ml}{2}(wz+xy)X_5+\frac{l}{2}\{1+m(x^2+z^2)\}X_6+\frac{ml}{2}(wx-yz)X_7,\\
\nabla_{X_2}X_5 & = \frac{ml}{2}(wz+xy)X_4-\frac{ml}{2}(wx-yz)X_6-\frac{l}{2}\{1+m(w^2+y^2)\}X_7,\\
\nabla_{X_2}X_6 & = -\frac{l}{2}\{1+m(x^2+z^2)\}X_4+\frac{ml}{2}(wx-yz)X_5+\frac{ml}{2}(wz+xy)X_7,\\
\nabla_{X_2}X_7 & = -\frac{ml}{2}(wx-yz)X_4+\frac{l}{2}\{1+m(w^2+y^2)\}X_5-\frac{ml}{2}(wz+xy)X_6,\\
\nabla_{X_3}X_4 & =\frac{ml}{2}(wy-xz)X_5-\frac{ml}{2}(wx+yz)X_6+\frac{l}{2}\{1+m(x^2+y^2)\}X_7,\\
\nabla_{X_3}X_5 & =-\frac{ml}{2}(wy-xz)X_4+\frac{l}{2}\{1+m(w^2+z^2)\}X_6-\frac{ml}{2}(wx+yz)X_7,\\
\nabla_{X_3}X_6 & =\frac{ml}{2}(wx+yz)X_4-\frac{l}{2}\{1+m(w^2+z^2)\}X_5-\frac{ml}{2}(wy-xz)X_7,\\
\nabla_{X_3}X_7 & =-\frac{l}{2}\{1+m(x^2+y^2)\}X_4+\frac{ml}{2}(wx+yz)X_5+\frac{ml}{2}(wy-xz)X_6,\\
\nabla_{X_4}X_4 & =2m(xX_5+yX_6+zX_7),\\
\nabla_{X_4}X_5 & =-\frac{l}{2}\{1+m(y^2+z^2)\}X_1+\frac{ml}{2}(wz+xy)X_2-\frac{ml}{2}(wy-xz)X_3-2mxX_4,\\
\nabla_{X_4}X_6 & =-\frac{ml}{2}(wz-xy)X_1-\frac{l}{2}\{1+m(x^2+z^2)\}X_2+\frac{ml}{2}(wx+yz)X_3-2myX_4,\\
\nabla_{X_4}X_7 & =\frac{ml}{2}(wy+xz)X_1-\frac{ml}{2}(wx-yz)X_2-\frac{l}{2}\{1+m(x^2+y^2)\}X_3-2mzX_4,\\
\nabla_{X_5}X_4 & =\frac{l}{2}\{1+m(y^2+z^2)\}X_1-\frac{ml}{2}(wz+xy)X_2+\frac{ml}{2}(wy-xz)X_3-2mwX_5,\\
\nabla_{X_5}X_5 & =2m(wX_4+yX_6+zX_7),\\
\nabla_{X_5}X_6 & =-\frac{ml}{2}(wy+xz)X_1+\frac{ml}{2}(wx-yz)X_2-\frac{l}{2}\{1+m(w^2+z^2)\}X_3-2myX_5,\\
\nabla_{X_5}X_7 & =-\frac{ml}{2}(wz-xy)X_1+\frac{l}{2}\{1+m(w^2+y^2)\}X_2+\frac{ml}{2}(wx+yz)X_3-2mzX_5,\\
\end{split}}
\end{equation*}

\begin{equation*}\footnotesize{
\begin{split}
\nabla_{X_6}X_4 & =\frac{ml}{2}(wz-xy)X_1+\frac{l}{2}\{1+m(x^2+z^2)\}X_2-\frac{ml}{2}(wx+yz)X_3-2mwX_6,\\
\nabla_{X_6}X_5 & =\frac{ml}{2}(wy+xz)X_1-\frac{ml}{2}(wx-yz)X_2+\frac{l}{2}\{1+m(w^2+z^2)\}X_3-2mxX_6,\\
\nabla_{X_6}X_6 & =2m(wX_4+xX_5+zX_7),\\
\nabla_{X_6}X_7 & =-\frac{l}{2}\{1+m(w^2+x^2)\}X_1-\frac{ml}{2}(wz+xy)X_2+\frac{ml}{2}(wy-xz)X_3-2mzX_6,\\
\nabla_{X_7}X_4 & =-\frac{ml}{2}(wy+xz)X_1+\frac{ml}{2}(wx-yz)X_2+\frac{l}{2}\{1+m(x^2+y^2)\}X_3-2mwX_7,\\
\nabla_{X_7}X_5 & =\frac{ml}{2}(wz-xy)X_1-\frac{l}{2}\{1+m(w^2+y^2)\}X_2-\frac{ml}{2}(wx+yz)X_3-2mxX_7,\\
\nabla_{X_7}X_6 & =\frac{l}{2}\{1+m(w^2+x^2)\}X_1+\frac{ml}{2}(wz+xy)X_2-\frac{ml}{2}(wy-xz)X_3-2myX_7,\\
\nabla_{X_7}X_7 & =2m(wX_4+xX_5+yX_6).
\end{split}}
\end{equation*}

\vspace{.5cm}

The curvature tensor:
\begin{equation*}
\begin{split}
R_{X_1X_4X_1X_4}=R_{X_1X_5X_1X_5} & =\frac{l^2 }{4}\{1+m(K+1)(y^2+z^2)\},\\
R_{X_1X_6X_1X_6}=R_{X_1X_7X_1X_7} & =\frac{l^2 }{4}\{1+m(K+1)(w^2+x^2)\},\\
R_{X_2X_4X_2X_4}=R_{X_2X_6X_2X_6} & =\frac{l^2 }{4}\{1+m(K+1)(x^2+z^2)\},\\
R_{X_2X_5X_2X_5}=R_{X_2X_7X_2X_7} & =\frac{l^2 }{4}\{1+m(K+1)(w^2+y^2)\},\\
R_{X_3X_4X_3X_4}=R_{X_3X_7X_3X_7} & =\frac{l^2 }{4}\{1+m(K+1)(x^2+y^2)\},\\
R_{X_3X_5X_3X_5}=R_{X_3X_6X_3X_6} & =\frac{l^2 }{4}\{1+m(K+1)(w^2+z^2)\},\\
R_{X_4X_5X_4X_5} & = 4m-3R_{X_1X_4X_1X_4} ,\\
R_{X_4X_6X_4X_6} & = 4m-3R_{X_2X_4X_2X_4} ,\\
R_{X_4X_7X_4X_7} & = 4m-3R_{X_3X_4X_3X_4} ,\\
R_{X_5X_6X_5X_6} & = 4m-3R_{X_3X_5X_3X_5} ,\\
R_{X_5X_7X_5X_7} & = 4m-3R_{X_2X_5X_2X_5} ,\\
R_{X_6X_7X_6X_7} & = 4m-3R_{X_1X_6X_1X_6} .\\
\end{split}
\end{equation*}

\vspace{1cm}

{\bf Acknowledgements.} We wish to thank Professors V. Miquel and F. J. Carreras for his valuable and enlightening comments on the subject of this paper. AF has been partially supported by MINECO/FEDER grant MTM2015-65430-P and Fundaci\'{o}n S\'{e}neca project 19901/GERM/15. AMN has been partially supported by MINECO-FEDER grant MTM2016-77093-P. ADT has been partially supported by Red IEMath-Galicia, reference CN 2012/077, Spain.


\end{document}